\documentclass[12pt]{amsart}
\usepackage{amssymb,amsmath}
\numberwithin{equation}{section}
\def\T{\text}
\def\1#1{\overline{#1}}
\def\2#1{\widetilde{#1}}
\def\3#1{\widehat{#1}}
\def\4#1{\mathbb{#1}}
\def\5#1{\frak{#1}}
\def\6#1{{\mathcal{#1}}}
\def\C{{\4C}}
\def\R{{\4R}}

\def\Z{{\4Z}}
\def\sumK{\underset{|K|=k-1}{{\sum}'}}
\def\simleq{\underset\sim<}
\def\sumJ{\underset{|J|=k}{{\sum}'}}
\def\sumij{\underset {ij=1,\dots,n-l}{{\sum}}}

\def\sumjq{\underset {j\leq q_o}\sum}

\def\sumj{\underset {j=1,\dots,n}{{\sum}}}

\begin{document}
\abstract
We prove non-subelliptic estimates for the tangential Cauchy-Riemann system over a weakly ``$q$-pseudoconvex" higher codimensional submanifold $M$ of $\C^n$. Let us point out that our hypotheses do not suffice to guarantee subelliptic estimates, in general. Even more: hypoellipticity of the tangential C-R system is not in question (as shows the example by {\sc Kohn} 
of \cite{K73} in case 
of a  Levi-flat hypersurface). However our estimates suffice for existence of smooth solutions to the inhomogeneous C-R equations in certain degree.

The main ingredients in our proofs are the weighted $L^2$ estimates by {\sc H\"ormander} 
\cite{H65} and {\sc Kohn} \cite{K73} of \S 2 and the tangential $\bar\partial$-Neumann operator 
by {\sc Kohn} of \S 4; for this latter we also refer to the book \cite{CS01}. 
As for the notion of $q$ pseudoconvexity we follow closely {\sc Zampieri} \cite{Z00}. 
The main technical result, Theorem~\ref{t2.1}, is a version for ``perturbed" $q$-pseudoconvex   domains of a similar result by {\sc Ahn} \cite{A03} who generalizes in turn {\sc Chen-Shaw} \cite{CS01}.  
\endabstract
\title[Non-subelliptic  estimates...]{Non-subelliptic  estimates for the tangential 
Cauchy-Riemann system} 
\author[H.~Ahn, L.~Baracco, G.~Zampieri]{Heungju Ahn, Luca Baracco, Giuseppe Zampieri}
\address{Dipartimento di Matematica, Universit\`a di Padova, via 
Belzoni 7, 35131 Padova, Italy}
\email{hjahn@math.unipd.it, baracco@math.unipd.it,  
zampieri@math.unipd.it}
%\subjclass{}
\maketitle
%\tableofcontents
% Standard sets
\def\Giialpha{\mathcal G^{i,i\alpha}}
\def\cn{{\C^n}}
\def\cnn{{\C^{n'}}}
\def\ocn{\2{\C^n}}
\def\ocnn{\2{\C^{n'}}}
% Abbreviations
\def\const{{\rm const}}
\def\rk{{\rm rank\,}}
\def\id{{\sf id}}
\def\aut{{\sf aut}}
\def\Aut{{\sf Aut}}
\def\CR{{\rm CR}}
\def\GL{{\sf GL}}
\def\Re{{\sf Re}\,}
\def\Im{{\sf Im}\,}
\def\codim{{\rm codim}}
\def\crd{\dim_{{\rm CR}}}
\def\crc{{\rm codim_{CR}}}
\def\phi{\varphi}
\def\eps{\varepsilon}
\def\d{\partial}
\def\a{\alpha}
\def\b{\beta}
\def\g{\gamma}
\def\G{\Gamma}
\def\D{\Delta}
\def\Om{\Omega}
\def\k{\kappa}
\def\l{\lambda}
\def\L{\Lambda}
\def\z{{\bar z}}
\def\w{{\bar w}}
\def\Z{{\1Z}}
\def\t{{\tau}}
\def\th{\theta}
\emergencystretch15pt
\frenchspacing
\newtheorem{Thm}{Theorem}[section]
\newtheorem{Cor}[Thm]{Corollary}
\newtheorem{Pro}[Thm]{Proposition}
\newtheorem{Lem}[Thm]{Lemma}
\theoremstyle{definition}\newtheorem{Def}[Thm]{Definition}
\theoremstyle{remark}                                                                                                                  
\newtheorem{Rem}[Thm]{Remark}
\newtheorem{Exa}[Thm]{Example}
\newtheorem{Exs}[Thm]{Examples}
\def\Label#1{\label{#1}}
\def\bl{\begin{Lem}}
\def\el{\end{Lem}}
\def\bp{\begin{Pro}}
\def\ep{\end{Pro}}
\def\bt{\begin{Thm}}
\def\et{\end{Thm}}
\def\bc{\begin{Cor}}
\def\ec{\end{Cor}}
\def\bd{\begin{Def}}
\def\ed{\end{Def}}
\def\br{\begin{Rem}}
\def\er{\end{Rem}}
\def\be{\begin{Exa}}
\def\ee{\end{Exa}}
\def\bpf{\begin{proof}}
\def\epf{\end{proof}}
\def\ben{\begin{enumerate}}
\def\een{\end{enumerate}}

\def\simto{\overset\sim\to\to}
\def\1alpha{[\frac1\alpha]}
\def\T{\text}
\def\R{{\Bbb R}}
\def\I{{\Bbb I}}
\def\C{{\Bbb C}}
\def\Z{{\Bbb Z}}
\def\Fialpha{{\mathcal F^{i,\alpha}}}
\def\Fiialpha{{\mathcal F^{i,i\alpha}}}
\def\Figamma{{\mathcal F^{i,\gamma}}}
\def\Real{\Re} 
%
%
%
%\endtopmatter
\section{  $q$-pseudoconvexity  in higher codimension}
Let $M$ be a real generic  submanifold of $\C^n$ of codimension $l$ 
and of class $C^i$, $i\geq4$, defined by a system of $l$ independent equations $\rho^h=0,\,\,h=n-l+1,\dots,n$.
We denote by $\rho$ the vector valued function with components $\rho^h$. Let $TM$ denote the tangent bundle to $M$, 
$T^\C M=TM\cap iTM$ the complex tangent bundle,
$T^{1,0}M$ and $T^{0,1}M$ the subbundles of $\C\otimes_\R T^\C M$ of holomorphic and antiholomorphic forms respectively. 
 Let $\mathcal L_\rho$, resp. $\mathcal L_M$, be the Levi form of 
$\rho$, resp. $M$, which is the Hermitian form defined, in a system of complex coordinates $z=x+iy$ for $\C^n$, by 
the matrix $\left(\partial^2_{z_i\bar z_j}\rho\right)$, resp. $\left.\left(\partial^2_{z_i\bar z_j}\rho
\right)\right|_{T^\C M}$. Let $T^*_M\C^n$ denote the conormal bundle to $M$ consisting of those $(0,1)$ forms whose real part vanish 
over $TM$ and set $\dot T^*_M\C^n= T^*_M\C^n\setminus\{0\}$. Note that the set of the $\partial\rho^h$'s are a basis for $T^*_M\C^n$. Identify in this basis  
$\R^l\mapsto \left(T^*_M\C^n\right)|_z\,\,z\in M$ by $a\mapsto \xi:=\sum_ha_h\partial\rho^h(z)$;
this yields an identification $T^*_M\C^n\simeq M\times \R^l$ and $\dot T^*_M\C^n/\R^+\simeq M\times S^{l-1}$ where $S^{l-1}$ is the spherical surface of dimension $l-1$.
 Set $\partial^2_{z_i\bar z_j}\rho^\xi=\sum_ha_h
\partial^2_{z_i\bar z_j}\rho^h$ and define $\mathcal L_\rho^\xi(z)=\left(\partial^2_{z_i\bar z_j}\rho^\xi(z)\right)_{ij}$ and
\begin{equation}
\Label{1.0}
\mathcal L_M^{\xi}(z)=\left.\left(\partial^2_{z_i\bar z_j}\rho^{\xi}(z)\right)_{ij}\right|_{T^\C M}\quad (z,\xi)\in M\times S^{l-1}.
\end{equation}
 The form $\mathcal L^\xi_M(z)$ is called the ``microlocal" Levi form of $M$ at $z$ in codirection $\xi$. Note that the Levi form is independent of the choice of a system of equations $\rho=0$ for $M$.  
 \bd
 \Label{d1.1}
We will deal with the assumption that there exists a smooth subbundle $\mathcal V^{q_o}=\mathcal V^{q_o}_{(x,\xi)}$ of $T^\C M$ of rank $q_o\leq q$ such that for any bundle $\mathcal V^{q+1}$ of rank $q+1$ we have
 \begin{equation}
 \Label{1.1}
 \T{trace}\left(\mathcal L_M^\xi(z)\right)\big|_{\mathcal V^{q+1}_{(z,\xi)}}-\T{trace}\left(\mathcal L_M^\xi(z)\right)
 \big|_{\mathcal V^{q_o}_{(z,\xi)}}\geq0\T{  $\forall (z,\xi)\in M\times S^{l-1}$}.
 \end{equation}
 \ed
We will deal also with the local version  of \eqref{1.1} at $z_o$ in which the condition  holds for any $(z,\xi)\in M'\times S^{l-1}$ where $M'$ is a neighborhood of $z_o$ in $M$. 
Let us denote by $\lambda_j=\lambda_j^\xi(z)$ the eigenvalues of $\mathcal L_M=\mathcal L_M^\xi(z)$
ordered as $\lambda_1\leq\lambda_2\leq\dots\leq\lambda_{n-l}$,
 and by $s^+=s^+(z,\xi),\,\,s^-=s^-(z,\xi),\,\,s^0=s^0(z,\xi)$ the numbers of its respectively positive, negative, and null eigenvalues; 
 note that $s^+(z,-\xi)=s^-(z,+\xi)$.
 We consider an orthonormal basis $\{\omega_j\}_{j\leq n}$ of $(1,0)$ forms and the dual basis $\{\partial_{\omega_j}\}_{j\leq n}$ of $(1,0)$ vector fields. We make our choice so that $\omega_j=\partial\rho^{j-n+l}$ for any $j\geq n-l+1$ and decompose the basis into $\{\omega'_j\}_{j\leq n-l}$ and $\{\omega''_j\}_{j\geq n-l+1}$, and use the similar decomposition for the dual basis $\{\partial'_{\omega_j}\}_{j\leq n-l}$ and  $\{\partial''_{\omega_j}\}_{j\geq n-l+1}$. If we change the partial basis $\{\omega'_j\}$ so that $\mathcal V^{q_o}=\T{Span}\{\Real \partial_{\omega_j}\}_{j\leq q_o}$,   then \eqref{1.1} reads as
\begin{equation}
\Label{1.2}
\sum_{j\leq q+1}\lambda_j^\xi(z)-\sum_{j\leq q_o}\rho_{jj}^\xi(z)\geq0\,\,\forall (z,\xi)\in M\times S^{l-1}.
\end{equation}
Note that  \eqref{1.2} implies the similar property with $q+1$ replaced by any $k\geq q+1$. In fact if \eqref{1.2} holds, then $\lambda_{q+1}\geq0$ and hence $\lambda_j\geq0\,\,\forall j\geq q+1$. Thus the terms $\lambda_j$ with $q+2\leq j\leq k$ can be added to the left hand side of the inequality without destroying it. In case $q=0$, condition \eqref{1.2} reduces to $\lambda_1\geq0$ which 
means that $\mathcal L_M$ is positive semi-definite: thus $M$ is pseudoconvex in the classical sense.
 For $\mathcal L_M$ non-degenerate and with $q=s^-$, we regain the classical notion of strong $q$-pseudoconvexity (cf. for instance \cite{H65}). In the weak case, that is when in \eqref{1.1} or \eqref{1.2} we have weak inequalities, our condition goes back to \cite{Ho92} and, more closely, to \cite{Z00}; it was also recently refined by \cite{A03}. Before entering algebraic details about \eqref{1.1} or \eqref{1.2} we wish to discuss some examples. In all of them we have $q=q_o$.
\be
\Label{e1}
We fix a point $z_o$ and set $q=\underset\xi{\T{sup}}(s^-(z_o,\xi)+s^0(z_o,\xi))$. Note that $s^-+s^0$ is upper semicontinuous in $z$ and $\xi$ and integer valued; hence it attains a local maximum at any point. In particular $q$ remains unchanged if we take the supremum also with respect to $z$ in a neighborhood of  $z_o$. We choose $\xi_o$ where $s^-+s^0$ attains the global maximum in $S^{l-1}$. Let $\mathcal V^q_{(z_o,\xi_o)}$ be the span of the negative and null eigenvectors at $(z_o,\xi_o)$ that we identify to the span of the first $q$ coordinate vectors. We have
\begin{equation}
\Label{a}
\sum_{j\leq q+1}\lambda_j^{\xi_o}(z_o)-\sum_{j\leq q}\rho_{jj}^{\xi_o}(z_o)\geq\lambda^{\xi_o}_{q+1}(z_o)>0.
\end{equation}
Hence if we move $(z,\xi)$ near $(z_o,\xi_o)$ it remains true, by continuity, that \eqref{a} is $>0$. In general, for any $(z,\xi)$, we can choose $\mathcal V^q_{(z,\xi)}$ such that $\sum_{j\leq q+1}\lambda_j^\xi(z)-\sum_{j\leq q}\rho_{jj}^\xi(z)>0$, though the difference in the left side needs not to coincide with $\lambda_{q+1}^\xi(z)$. By a partition argument over  the unit circle $S^{l-1}$, we get \eqref{a} for all $(z,\xi)\in M'\times S^{l-1}$
where $M'$ is a neighbourhood of $z_o$. 
 The above condition is considered by {\sc Nacinovich} in \cite{N84} where existence theorems for the tangential $\bar\partial$ system are derived.
Our task is to refine the above criterion and move $q$ to lower values.
\ee
\be
\Label{e2}
Let $M$ be a hypersurface; in this situation $M\times S^0$ consists of just two components $(z,\pm\xi)$. We write $\lambda_j^\pm(z)$ instead of $\lambda_j^{\pm\xi}(z)$, $s^\pm(z)$ instead of $s^\pm(z,+\xi)$ and so on; note that $s^\pm(z)=s^-(z,\mp\xi)$. In this situation \eqref{1.2} means the existence of two bundles $\mathcal V^{q_o}_+$ and  $\mathcal V^{q_o}_-$ resp., such that in the two systems in which these bundles are reduced to the span of the first $q_o$ coordinate vectors, we have 
\begin{equation}
\Label{b}
\sum_{j\leq q+1}\lambda_j^\pm(z)-\sum_{j\leq q_o}\rho_{jj}^\pm(z)\geq0.
\end{equation}
According to Example~\ref{e1} a first rough index $q$ for which \eqref{1.2} holds in $M'\times S^{l-1}$ for a neighbourhood $M'$ of $z_o$, is
\begin{equation}
\Label{c}
q=\T{sup}(s^-(z_o),s^+(z_o))+s^0(z_o).
\end{equation}
In some cases we can do better. For instance, assume that $s^-(z)$ is constant for $z$ close to $z_o$. Then $\lambda^+_{s^-}<0\leq \lambda^+_{s^-+1}$ and hence the negative eigenvectors span a bundle $\mathcal V^{s^-}_+$ that,
identified to the span of the first $s^-$ coordinate vectors, yields $\sum_{j\leq s^-+1}\lambda^+_j(z)-\sum_{j\leq s^-}\rho_{jj}^+(z)\geq0$.
 Of course, the same can be said in case $s^+(z)$ is constant. 
 For the bundle $\mathcal V^{s^+}_-$ of the positive eigenvectors, identified to the first $s^+$ vectors, we have $\sum_{j\leq s^++1}\lambda^-_j(z)-\sum_{j\leq s^+}\rho^-_{jj}(z)\geq0$.
 Thus if both $s^\pm(z)$ are constant, or equivalently if the corank $s^0(z)$ is constant, then we have \eqref{1.2} at $z_o$ for 
$$
q=\T{sup}(s^-,s^+).
$$
Thus we succeeded in decresing by $s^0$ the value of $q$ with respect to \eqref{c}.
\ee 
We want to consider a variant of conditions \eqref{1.1} or \eqref{1.2} that we need first to express in new terms.
For ordered multiindices $J=j_1<j_2<...<j_k$ of length $|J|=k$, let us consider vectors $w=(w_J)$. Decompose $J=jK$ with $|K|=k-1$ and write $w_{jK}=\T{sign}{{jK}\choose J}w_J$ where ${{jK}\choose J}$ is the permutation which orders $jK$. 
%$w'_{\cdot K}=(w_{iK})_{i\leq n-l},\,\,w''_{\cdot K}=(w_{iK})_{i\geq n-l+1}$. 
We will deal with the class of tangential forms; these are the forms $w=(w_J)$ such that any coefficient $w_J$ is $0$ if $J$ contains some index $j=n-l+1,\dots,n$. Sometimes we denote these forms by the alternative notation $w^{\tau}$.
We will denote by ${{\sum}'}$ summation over ordered indices. 
One checks that  \eqref{1.2} is equivalent to
\begin{multline}
\Label{1.2bis}
\sumK{\underset {ij\leq n-l}\sum}\rho_{ij}^\xi(z)w_{iK}\bar w_{jK}-\sumJ\underset{j\leq q_o}\sum\rho_{jj}^\xi|w_J|^2\geq0
\\
\T{for any tangential form $w$ of length $k\geq q+1$, and $\forall (z,\xi)\in M\times S^{l-1}$}.                              
\end{multline}
Along with \eqref{1.2bis} we will also consider the condition
\begin{multline}
\Label{1.2ter}
\sumK{\underset {ij\leq n-l}\sum}\rho_{ij}^\xi(z)w_{iK}\bar w_{jK}-\sumK\underset{j\leq q_o}\sum\rho_{jj}^\xi(z)|w_{jK}|^2\geq0
\\
\T{for any tangential form $w$ of length $k\geq q+1$, and $\forall (z,\xi)\in M\times S^{l-1}$}. 
\end{multline}
One can also consider some intermediate condition between \eqref{1.2bis} and \eqref{1.2ter} in which for part of the indices $j\leq q_o$ one takes $J=jK$ and for the remaining indices one takes all $J$ without requiring $j\in J$. 
For $q=\underset\xi{\T{sup}}\,\,(s^-+s^0)$, 
\eqref{1.2bis} holds according to Example~\ref{e1}; in this case one sees that  \eqref{1.2ter} is also fulfilled.
 But we can also discuss some cases of \eqref{1.2ter} which do not fit \eqref{1.2bis}.
\be
\Label{e3}
Let $M_1\times M_2\subset \C^{n_1}\times\C^{n_2}$ be quadric hypersurfaces  given by diagonal equations. Thus $\mathcal L^\pm_{M_i}\,\,i=1,2$ are diagonal at any $z$. We define $q_i=\T{sup}(s^+_{M_i},s^-_{M_i})$, denote by 
$\mathcal V^{q_i}$ the span in $T^{1,0}M$ of the non-null eigenvectors, and put
$$
q=\T{sup}(n_1-1+q_2,n_2-1+q_1).
$$
Consider a point, say $\xi=(0,+\xi_2)$, and take $\mathcal V^{n_1+q_2}=T^\C M_1\oplus\mathcal V^{q_2}$
where $\mathcal V^{q_2}$ contains the span of the $s^-$ negative eigenvectors. 
 We have
\begin{multline}
\sumK\,\,{\underset {ij\leq n_1-1}\sum}\rho_{ij}^\pm(z)w_{iK}\bar w_{jK}-\sumK\sum_{j\leq n_1-1}\rho^\pm_{jj}(z)|w_{jK}|^2
\\
+\sumK\,\,{\underset { n_1+1\leq i\,j\leq n_1+n_2-2}\sum}\rho_{ij}^+(z)w_{iK}\bar w_{jK}-\sumK\sum_{n_1+1\leq j\leq n_1+q_2}\rho^+_{jj}(z)|w_{jK}|^2
\\
\geq\lambda^+_{q_2+1}|w_{n_1+q_2+1}|^2+\dots+\lambda^+_{k-n_1}|w_k|^2\geq0.
\end{multline}
The above discussion applies for instance to the manifold in $\C^{n_1+n_2}$ defined by the equations
\begin{equation}
\begin{cases}
y_{n_1}=|z_1|^2-|z_2|^2,
\\
y_{n_1+n_2}=|z_{n_1+1}|^2-|z_{n_1+2}|^2.
\end{cases}
\end{equation}
Here $q_i=1$ for $i=1,2$ and therefore \eqref{1.2ter} is satisfied for $q=sup(n_1,n_2)+1$.
\ee

We refine now our choice of the basis of $(1,0)$ forms to make it better adapted to $M$. We first choose our equations $\rho^h=0\,\,h=n-l+1,...,n$ having orthonormal differentials along $M$. We extend the system of the $\partial\rho^h|_M$ to an orthonormal system $\{\omega''\}$ which spans $\T{Span}\{\partial\rho^h\}$ even ouside of $M$. We then take an orthonormal completion $\{\omega'\}$ of $\{\omega''\}$ and denote by $\{\partial'_\omega,\,\,\partial''_\omega\}$ the dual system of $(1,0)$ vector fields. Note that by our choice we have $\forall j,k$
\begin{equation} 
\Label{1.3}
\partial'_{\omega_j}\rho^k\equiv0\,\,\partial''_{\omega_j}\rho^k\equiv\varkappa_{jk}\,\,
\T{on $\C^n$},
\end{equation}
where $\varkappa_{jk}$ is  the Kronecker symbol.
We introduce now the spaces of forms of type $(0,k)$; in a basis $\{\omega_j\}$ they can be 
written as $u=\sumJ u_J\bar\omega_J\,\,|J|=k$ with coefficients in spaces of various kind 
such as $C^\infty(\C^n)$ or $L^2(\C^n)$ or, for a positive function $\phi$, $L^2_\phi(\C^n)$ 
that is the space of functions which satisfy $||u_J||_\phi:=\left(\int e^{-\phi}|u_J|^2dV
\right)^{\frac12}<\infty$. We will denote by $C^{\infty}_k,\,\,L^{2}_k,\,\,(L^{2}_{\phi})_k$ 
the above defined spaces,
and also denote them by the common symbol $\Lambda_k$ when we want not to stress attention to 
the kind of the coefficients.
We denote by $\mathcal C_k$ the 
restriction to $M$ of the
ideal of $\Lambda_k$ engendred by $\rho$ and $\bar\partial\rho$, 
and define the space of tangential forms on $M$ 
as the orthogonal complement of $\mathcal C_k$ in $\Lambda_k$:
$$
\mathcal T_k=\mathcal C_k^\perp. 
$$                                                               
Any tangential form can be represented as the restriction to $M$ of a form 
satisfying the $\bar\partial$-Neumann conditions on $M$:
$$
\sumj\left.\rho^h_ju_{jK}\right|_M=0\quad\forall h=n-l+1,...,n\,\,\forall K,
$$
where we have used the notation $\rho^h_j$ for $\partial_{\omega_j}\rho^h$. 
Let us take an orthonormal frame $\{\omega_j\}$ satisfying the above conditions and in particular 
\eqref{1.3}.
Let us decompose any $u$ as $u=u^{\tau}+u^\nu$ where in $u^\tau$ we collect coefficients 
corresponding to indices $J$ such that  $n-l+1,\dots,n\notin J$ and in $u^\nu$ the remaining 
ones. The fact that $u|_M$ satisfies the $\bar\partial$-Neumann conditions reads 
\begin{equation}
\Label{1.3bis}
u^\nu(z)|_M\equiv0\quad\forall z\in M.
\end{equation}
It is obvious that \eqref{1.3bis} implies 
$\partial'_{\omega_j}u^\nu|_M\equiv0\,\,\forall K$. 
We can see that we may  choose, among representatives of a tangential form, one which 
satisfies
\begin{equation}
\Label{1.3ter}
u^\nu(z)\equiv0\quad\forall z\in\C^n.
\end{equation}
We also choose the extension of $u$ from $M$ to $\C^n$ so that for all coefficients we have 
that $\partial''_{\bar\omega_j}u_J|_M$ is nearly $0$ according to the following 
considerations.
\bp
\Label{p1.0,5}
Let $M$ have class $C^{m+1}$ and be locally defined at $0$ by $y''_h=g_h$ for $g_h(0)=0,\,\,
\partial g_h(0)=0$. Then there is a 
local system of vector fields $\bar L_h\,\,h\geq n-l+1$  of class $C^m$ and type $(0,1)$ 
with $\bar L_h(0)=\partial_{\bar z_h}$ such that for any function 
$f\in C^m(M)$ there exists an extension $\tilde f$ in $\C^n$ such that
\begin{equation*}
\begin{cases}
\bar L_h\tilde f|_M\equiv\mathcal O^m,
\\
\tilde f|_M\equiv f,
\end{cases}
\end{equation*}
(where the symbol $\mathcal O^m$ denotes an infinitesimal of order $m$ with respect to the 
distance to $M$). 
\ep
\bpf
We consider the parametrization of $M$:
\begin{equation*}
G: \C^{n-l}\times \R^l\to M,\,\, (z',x'')\mapsto (z',x''+ig(z',x'')).
\end{equation*}
We extend $G$ to $\tilde G$ which is {\it m-holomorphic} along $M$ that is
\begin{equation*}
\tilde G: \C^{n-l}\times\C^l\to \C^n,\,\,(z',z'')\mapsto(z',z''+i\tilde g(z',z'')),
\end{equation*}
such that
\begin{equation}
\Label{Whitney1}
\partial_{\bar z''_h}\tilde g|_{\C^{n-l}\times\R^l}=\mathcal O^m.
\end{equation}
This statement belongs to the family of Whitney's extension 
theorems. Given $f$, we define $f_o:=f\circ G$, extend $f_o$ to $\tilde f_o$ from 
$\C^{n-l}\times\R^l$ to $\C^{n-l}\times\C^l$ with the property
\begin{equation}
\Label{Whitney2}
\partial''_{\bar z}\tilde f_o|_{\C^{n-l}\times\R^l}=\mathcal O^m,
\end{equation}
and set $\tilde f:=\tilde f_o\circ\tilde G^{-1}$. We also define
$$
\bar L_h=\tilde G_*\partial''_{\bar z_h}\quad h\geq n-l+1.
$$ 
It is clear that each $\bar L_h$ is of type $(0,1)$ along $M$ due to \eqref{Whitney1}. 
We also have
\begin{equation}
\begin{cases}
\bar L_h\sim\partial''_{\bar z_h}\T{ due to $\partial g_j(0)=0\,\forall j$},
\\
\bar L_h\tilde f|_M=\partial''_{\bar z_h}\tilde f_o|_{\C^{n-l}\times \R^l}=\mathcal O^m
\T{ due to \eqref{Whitney2}}.
\end{cases}
\end{equation}
\epf
\br
\Label{r1.1}
Let $\{\omega_j\}$ be an orthonormal system of $(1,0)$-forms with $\omega''_h=\partial\rho^h,\,\,h\geq n-l+1$ for a system of equations $\rho^h=0$ of $M$ such that $\partial\rho^h(0)=dz_h$.
If $z\in\C^n$ is close to $M$ and $z^*$ is the point on $M$ of minimal distance, we have
\begin{equation*}
\partial_{\bar\omega_k(z)}=\partial_{\bar \omega_k(z^*)}+O(|z-z^*|), \,\,\forall k \T{ such that }1\leq  
k\leq n.
\end{equation*}
For $h\geq n-l+1$, write 
$$
\partial''_{\bar\omega_h}|_M=\sum_{i\geq n-l+1}b_i\bar L_i|_M+\sum_{j\leq n-l}a_j
\partial'_{\bar\omega_j}|_M;
$$
note that the $a_j$'s and the $b_i$'s for $i\neq h$ are small. Thus, if $\tilde f$ satisfies 
the conclusions of the preceding proposition we have
\begin{equation}
\Label{normal}
\partial''_{\bar\omega_h}\tilde f(z)=\sum_{j\leq n-l}a_j\partial'_{\bar\omega_j}
f(z)+O(|z-z^*|),
\end{equation}
for small $a_j$'s. 
\er

 We let
$$
{\tilde\rho}=\frac{|\rho|^2-\epsilon^2}{2\epsilon},
$$
and  define the system of ``tuboidal" neighbborhoods of $M$ adapted to the frame $\omega$ by
$$
U_\epsilon=\{z\in\C^n:\,{\tilde\rho}(z)<0\}.
$$
Let $|\rho(z)|=\epsilon$ and $a=\epsilon^{-1}(\rho^h(z))_h$; recall that we are identifying $a$  to a cotangent vector $\xi$  the one with coordinates $a$ in the system of $1$-forms $\partial\rho^h,\,\,h=n-l+1,...,n$; note that $\xi$ is conormal to $\partial U_\epsilon$. Let $\C^n\to M,\,\,z\mapsto z^*$ be any transversal projection. We have, if $\rho$ is $C^k$ and keeping the assumption that $z$ belongs to $\partial U_\epsilon$ and $z^*$ is the point of minimal distance on $M$ all through the sequel:
\begin{equation}
\Label{1.4}
|\rho^h_{ij}(z)-\rho^h_{ij}(z^*)|=O(\epsilon)\,\,\forall h,\,\forall ij.
\end{equation}
We also have
\begin{equation}
\Label{1.5}
\begin{split}
\mathcal L_{{\tilde\rho}}&=\epsilon^{-1}\sum_h\partial\rho^h\otimes\bar\partial\rho^h+\epsilon^{-1}\sum_h\rho^h\mathcal L_{\rho^h}
\\
&=\epsilon^{-1}\sum_h\partial\rho^h\otimes\bar\partial\rho^h+\mathcal L_\rho^\xi.
\end{split}
\end{equation}
We write $u'_{\cdot K}=(u_{iK})_{i\leq n-l},\,\,u''_{\cdot K}=(u_{iK})_{i\geq n-l+1}$. It follows for any $K$
\begin{equation}
\Label{1.6}
\mathcal L_{{\tilde\rho}}(u_{\cdot K},\bar u_{\cdot K})\geq \epsilon^{-1}|u''_{\cdot K}|^2+\mathcal L_\rho^\xi(u'_{\cdot K},\bar u'_{\cdot K})-c_1|u'_{\cdot K}||u''_{\cdot K}|-c_2|u''_{\cdot K}|^2,
\end{equation}
and hence
\begin{equation}
\Label{1.7}
\mathcal L_{\tilde\rho}(z)(u_{\cdot K},\bar u_{\cdot K})\geq \frac{\epsilon^{-1}}2 |u''_{\cdot K}|^2+\mathcal L_\rho^\xi(z^*)(u'_{\cdot K},\bar u'_{\cdot K})-O(\epsilon)|u'_{\cdot K}|^2.
\end{equation}
By combining \eqref{1.4} and \eqref{1.7} and by taking summation on $K$, we get the proof of the following statement which describes how \eqref{1.2bis} is affected when $z$ is no more a point of $M$, and $u$ is not necessarily a tangential form.  
\bt
\Label{t1.1}
Let $M$ satisfy \eqref{1.2}; then
\begin{multline}
\Label{new}
\sumK\sumij\tilde\rho_{ij}(z) u_{iK}\bar u_{jK}-\sumJ\sumjq\tilde\rho_{jj}(z)|u_J|^2
\geq-O(\epsilon)|u|^2
\\
\forall z\in\partial U_\epsilon\T{ and $\forall u$ of length $k\geq q+1$}.
\end{multline}
\et
One has also a local version of Theorem~\ref{t1.1} in a neighborhood of $z_o$.
\br
\Label{r1.0} 
If, instead of \eqref{1.2}, we assume \eqref{1.2ter}, then we have the similar conclusion as \eqref{new} but with the second term in the left containing only the indices for which $j\in J$, or equivalently those  in the form $J=jK$.
\er 
\br
\Label{r1.2}
The coefficients $a_{ij}$ of the basis of forms $\{\omega_i\}$ in which \eqref{new} holds are singular in $M$. In particular forthe normal vector fields we have that $(\partial''_{\omega_h}+\partial''_{\bar\omega_h})a_{ij} \,\forall h\geq n-l$ grow as $|\rho|^{-1}$.  However, for the tangential vector fields, we have that  $\partial'_{\omega_k}a_{ij},\,\,\partial'_{\bar\omega_k}a_{ij}\,\,\forall k\leq n-l$ and  $\frac1{2i}(\partial_{\omega_h}-\partial_{\bar\omega_h})a_{ij}\,\forall h\geq n-l+1$ are bounded. 
\er

\br
It follows from Theorem~\ref{1.1} that $M$ has a fundamental system of neighborhoods which are ``almost" $q$-pseudoconvex. In general these neighborhoods cannot be $q$-pseudoconvex as shows the example by {\sc Diederich-Fornaess} of non-trivial ``nebenh\"ulle".
\er 

Recall that $\partial'_{\omega_k}\rho^h\equiv0$ and that $\partial''_{\omega_k}\rho^h=\varkappa_{hk}$.
 It follows
\begin{equation}
\begin{split}
\sum_h\rho_{ij}^h{\tilde\rho}_h|_z&=\sum_h\rho^h_{ij}\epsilon^{-1}\sum_k\rho^k\partial_{\omega_h}\rho^k|_z
\\
&=\sum_h\rho^h_{ij}\epsilon^{-1}\rho^h|_z=\rho^\xi_{ij}(z)\,\,\T{for $z\in\partial U_\epsilon$} \T{ and $\xi:=\epsilon^{-1}\sum_h\rho^h\partial\rho^h$}.
\end{split}
\end{equation}
 Choose any transversal projection $z\mapsto z^*$; we have
\begin{equation}
\Label{1.8}
\sum_h\rho_{ij}^h(z)\tilde\rho_h(z)(u,\bar u)=\rho_{ij}^\xi(z^*)(u,\bar u)+O(\epsilon)|u|^2\,\,\T{for $z\in \partial U_\epsilon$}.
\end{equation}
This gives the proof of the following
\bp
\Label{p1.1}
Let $M$ satisfy \eqref{1.2}. Then 
\begin{multline}
\Label{newnew}
\sumK\sumij\sum_h\rho^h_{ij}(z)\tilde\rho_h(z) u_{iK}\bar u_{iK} \\
-\sumJ\sumjq\sum_h\rho_{jj}(z){\tilde\rho}_h|u_J|^2\geq-O(\epsilon)|u|^2,
\end{multline}
for any $ z\in\partial U_\epsilon$  and
for any form $u=u(z)$ (not necessarily satisfying $\bar\partial$-Neumann conditions on $M$) of order $k\geq q+1$.
\ep
Again, we have a local version at $z_o$ of this statement and also a variant under the assumption \eqref{1.2ter}.
\bd
\Label{defqpseudo}
We say that $M$ is $q$-pseudoconvex, resp. locally $q$-pseudoconvex at $z_o$,   when \eqref{1.1} or \eqref{1.2ter} are fulfilled for any $(z,\xi)\in M\times S^{l-1}$, resp. for any $(z,\xi)\in M'\times S^{l-1}$ for a neighborhood $M'$ of $z_o$.
\ed
\section{$L^2$ estimates for the ambient $\bar\partial$ system}
We denote by $u(z)=(u_J(z))\quad z\in M\subset \C^n$, a  form of type $(0,k)$ satisfying the $\bar\partial$-Neumann conditions; most of times its coefficients are supposed to be smooth. We also suppose that the orthonormal frame $\{\omega',\omega''\}$ and the extension $u$ satisfy all conditions listed in \S 1 including Proposition~\ref{p1.0,5} and the related remark. In particular $u^\nu\equiv0$ also outside $M$ and 
\begin{equation}
\Label{2.0}
\partial''_{\bar\omega_h}u_J=\sum_{j\leq n-l}a_j\partial'_{\bar\omega_j}u_J+O(|\rho|),
\end{equation}
 with small coefficients $a_j$.   
We denote by $||\cdot||_{H^0(M)}$ or $||\cdot||_{H^0(U_\epsilon)}$ the $H^0=L^2$ norms on $M$ and $U_\epsilon$ respectively; for any real positive function $\phi$ we denote by $H^0_\phi$ the $L^2$ norms with weight $e^{-\phi}$. We will make our choice of $\phi$ as $\phi=(t+c)|z|^2$ for a large parameter $t$ and for a constant $c$ depending on the coefficients of the $\omega_j$'s. 
 We denote by $\bar\partial$ , resp. $\bar\partial'$, the complex on antiholomorphic forms associated to all antiholomorphic vector fields $\partial_{\bar\omega_j},\,\,1\leq j\leq n$, resp. to $\partial'_{\bar\omega_j},\,\,1\leq j\leq n-l$. We denote by $\bar\partial^*$, resp. $\bar\partial^{\prime *}$ the $H^0_\phi$-transposed; note that $\bar\partial^*=\bar\partial^{\prime*}+O(|\rho|)$ over $\bar\partial$-Neumann forms. 
 We will still denote by $U_\epsilon$ the intersection of the tube $U_\epsilon$ with a suitable sphere centered at $z_o$.
 \bt
\Label{t2.1}
Let $M$ be $q$-pseudoconvex at $z_o$. Then for any $\bar\partial$-Neumann form $u$ of degree $k\geq q+1$ with support   whose coefficients satisfy \eqref{2.0}, and for any large real $t$, we have
\begin{equation}
\Label{2.1}
\frac t2||u||^2_{H^0_\phi(U_\epsilon)}\leq ||\bar\partial^{\prime*} u||^2_{H^0_\phi(U_\epsilon)}+||\bar\partial' u||^2_{H^0_\phi(U_\epsilon)}+o(\epsilon^l).
\end{equation}
\et
\bpf
We will only give the proof under the asumption  \eqref{1.2bis} in local form, the case of \eqref{1.2ter} being analogous. We also point out that by cutting the tube $U_\epsilon$ by a sphere we still have a domain which satisfies \eqref{new} and \eqref{newnew} in each smooth part of the boundary. Also, in the integrations by parts, some integrals in the 2-codimensional strata appear. But these are positive and so we can neglect them or equivalently we can assume from the beginning that $U_\epsilon$ is compact, smooth and satisfies \eqref{new} and \eqref{newnew}. 
The proof is closely related to that by Ahn \cite{A03} who deals with a $q$ pseudoconvex domain and gets the similar estimate as \eqref{2.1} without the error term $o(\epsilon^l)$. 
We simplify our notation and write $||\cdot||_\phi$ instead of $||\cdot||_{H^0_\phi(U_\epsilon)}$ all through the proof.
We also drop the symbol ${}^\prime$ in most of  notations: it will be understood that our indices will generally vary between $1$ and $n-l$.
We set $\phi_j=\partial_{\omega_j}\phi$ and define
$$
\delta_{\omega_j}=\partial_{\omega_j}-\phi_j;
$$
Hence $\delta_{\omega_j}$ is the transposed of $-\partial_{\bar\omega_j}$ in the weighted $H^0_\phi$ scalar product apart from a $0$-order operator which depends on tangential derivatives of the coefficients of the forms $\omega_j$'s.
We have
\begin{multline}
\Label{2.2}
\sumK\sumij\int_{U_\epsilon}e^{-\phi}\left(\delta_{\omega_i}u_{iK}\overline{\delta_{\omega_j}u_{jK}}-\partial_{\bar\omega_j}u_{iK}\overline{\partial_{\bar\omega_i}u_{jK}}\right)dV
\\
+\sumJ\sumj\int_{U_\epsilon}e^{-\phi}|\partial_{\bar\omega_j}u_J|^2 dV   \leq2(||\bar\partial^{\prime*}u||^2_\phi+||\bar\partial' u||^2_\phi)+R^1,
\end{multline}
where $R^1$ is an error term which only involves integration of $|u|^2$ and not of its derivatives. We will use the notation ``s.c.", resp. ``l.c.", to denote small constants, resp. large constants. 
We  have
\begin{equation}
\Label{2.3}
||\partial'_{\partial\omega_j}u_J||^2_\phi=||\delta_{\omega_j}'u_J||^2_\phi+\int_{U_\epsilon}e^{-\phi}\left[\delta_{\omega_j}',\partial_{\bar\omega_j}\right]u_J\bar u_J dV\,+\,R^2_{jj}\quad\forall j\leq n-l,
\end{equation}
where $R^2_{jj}$ can be estimated  both by   $\T{s.c.}||\partial_{\bar\omega_j}u_J||^2_\phi+\T{l.c.}||u_J||^2_\phi$ or $\T{s.c}||\delta_{\omega_j}u_J||^2_\phi+\T{l.c.}||u||^2_\phi$. In fact  the boundary integrals which arise in the integrations by parts for interchanging $\partial'_{\partial\omega_j}$ with $\delta'_{\omega_j}$ are $0$ due to  \eqref{1.3} that is $\partial'_{\omega_j}\rho\equiv0$. 
 We rewrite now the 
integrals of $\partial_{\bar\omega_j}u_{iK}\overline{\partial_{\bar\omega_i}u_{jK}}$ in the left side of \eqref{2.2}.
Integration by parts yields
\begin{multline}
\Label{2.4}
\int_{U_\epsilon}e^{-\phi}\partial_{\bar\omega_j}u_{iK}\overline{\partial_{\bar\omega_i}u_{jK}}=\int_{+\partial U_\epsilon}e^{-\phi}\partial_{\bar\omega_j}(u_{iK})\tilde\rho_i\bar u_{jK} dV
\\
-\int_{U_\epsilon}e^{-\phi}\delta_{\omega_i}\partial_{\bar\omega_j}(u_{iK})\bar u_{jK} dV\,+\,R^3_{ij},
\end{multline}
where $R^3_{ij}$ is an error which involves integrals of $\bar u_{jK}\partial_{\bar\omega_j}u_{iK}$. Again, in \eqref{2.4} the boundary integral is $0$: in fact, since $i\leq n-l$, then ${\tilde\rho}_i\equiv0$ (where we are using 
as always
the notation ${\tilde\rho}_i=\partial_{\omega_i}{\tilde\rho}$).
We also have
\begin{multline}
\Label{2.5}
\int_{U_\epsilon} e^{-\phi}\delta_{\omega_i}u_{iK}\overline{\delta_{\omega_j}u_{jK}} dV=\int_{+\partial U_\epsilon}e^{-\phi}\bar{\tilde\rho}_j\bar u_{jK}\delta_{\omega_i}u_{iK} dV
\\
-\,\int_{U_\epsilon}e^{-\phi}\partial_{\bar\omega_j}\delta_{\omega_i}u_{iK}\bar u_{jK} dV\,+\,R^4_{ij},
\end{multline}
where $R^4_{ij}$ involves integrals of   $\delta_{\omega_i}u_{iK}\bar u_{jK}$. Again, the boundary integral in \eqref{2.5} is $0$ due to \eqref{1.3} and \eqref{1.3ter}. Thus in the left side of \eqref{2.2} we use \eqref{2.4}, \eqref{2.5} in the first two terms for any $i$ and $j$ and next \eqref{2.3} in the third, but now only for  $j\leq p$. In this way we can rewrite the left side of \eqref{2.2} as
\begin{multline}
\Label{2.6}
\left(\sumK\sumij\int_{U_\epsilon}e^{-\phi}\left[\delta_{\omega_i}, 
\partial_{\bar\omega_j}\right]u_{iK}\bar u_{jK} dV  \right.\\
\left.-\sumJ\sumjq\int_{U_\epsilon}e^{-\phi}\left[\delta_{\omega_j},\partial_{\bar\omega_j}\right]u_J\bar u_J dV \right)
\end{multline}
\begin{equation*}
+\left(\sumJ\sumjq\int_{U_\epsilon}e^{-\phi}|\delta_{\omega_j}u_J|^2 dV+\sumJ\sum_{j\geq p+1}\int_{U_\epsilon}e^{-\phi}|\partial _{\partial\omega_j}u_J|^2 dV\right)+R^5,
\end{equation*}
where $R^5$ is the sum of the $R^2_{jj}$'s (for $j\leq q$), the $R^3_{ij}$'s and the 
$R^4_{ij}$'s. We denote by $S$ the second term in \eqref{2.6} that is $(\sumJ\sumjq\cdot+
\sumJ\underset{j\geq p+1}\sum\cdot)$. The terms $R^{2}$ were already estimated. As for the 
remaining we clearly have an analogous estimate
\begin{equation}
\Label{2.7}
R^i\leq \T{s.c.}S+\T{l.c.}||u||^2_\phi \,\,\forall i\geq2.
\end{equation}
Clearly an estimate of the type \eqref{2.7} also holds for $R^1$. We pass now to compute the 
commutators $\left[\delta_{\omega_j},\partial_{\bar\omega_j}\right]$. 
Let $(c_{ij}^h)$ be the matrix of the 2-form $\partial_{\bar\omega_h}$; note that since for 
$h\geq n-l+1$, we have $\omega_h=\partial\rho^h$, then $(c_{ji}^h)=(\rho^h_{ij})$ 
is the matrix of the Levi-form $\mathcal L_{\rho^h}$. The identity $\bar\partial\partial=
-\partial\bar\partial$ yields
\begin{equation}
\left[\partial_{\omega_i},\partial_{\bar\omega_j}\right]=\sum_{h=1}^nc_{ji}^h
\partial_{\omega_h}-\sum_{h=1}^n\bar c_{ij}^h\partial_{\bar\omega_h}.
\end{equation}
We denote by $(\phi_{ij})$ the matrix of $\mathcal L_\phi$ which coincides, 
up to an error term, with $(t+2c)\varkappa_{ij}$. We get
\begin{equation}
\Label{commutator}
\begin{split}
\left[\delta_{\omega_i},\partial_{\bar\omega_j}\right]&=\left[\partial_{\omega_i}-\phi_i,
\partial_{\bar\omega_j}\right]
\\
&=\left[\partial_{\omega_i},\partial_{\bar\omega_j}\right]-\left[\phi_i,\partial_{\bar\omega_j}
\right]
\\
&=\phi_{ij}+\sum_hc_{ji}^h\delta_{\omega_h}-\sum_h\bar c_{ij}^h\partial_{\bar\omega_h}
\\
&=\phi_{ij}+\sum_{h\geq n-l+1}\rho^h_{ij}(\delta_{\omega_h}-\partial_{\bar\omega_h})+
\sum_{h\leq n-l}(c^h_{ji}\delta_{\omega_h}-\bar c_{ij}^h\partial_{\bar\omega_h}).
\end{split}
\end{equation}

Integration by parts yields, on account of \eqref{1.3}:
$$
\int_{U_\epsilon}e^{-\phi}|\bar c_{ij}^h\partial_{\bar\omega_h}u_J\bar u_I| dV\leq\T{l.c.}||u||^2_\phi+\T{s.c}S\,\,\forall h\leq n-l,
$$
and
$$
\int_{U_\epsilon}e^{-\phi}|c_{ji}^h\delta_{\omega_h}u_J\bar u_I| dV\leq\T{l.c.}||u||^2_\phi+\T{s.c.}S\,\,\forall h\leq n-l.
$$
For $h\geq n-l+1$, we want to interchange  $\delta_{\omega_h}$  with 
$\partial_{\bar\omega_h}$ in our integrals; we have
\begin{multline}
\Label{2.10}
\underset{h\geq n-l+1}\sum\int_{U_\epsilon}e^{-\phi}\rho_{ij}^h\delta_{\omega_h}u_J\bar u_I dV=
\underset{h\geq n-l+1}\sum\int_{+\partial U_\epsilon} e^{-\phi}\rho_{ij}^h\rho_h^\epsilon u_J
\bar u_I dV
\\
- \underset{h\geq n-l+1}\sum\int_{U_\epsilon}e^{-\phi}\rho_{ij}^hu_J\overline{\partial_{\bar\omega_h} u_I} dV  +R^7,
\end{multline}
Here, for the error term we have the estimate $R^7\leq c||u||^2_\phi+o(\epsilon^l)$. In fact the coefficients of the vector fields $\partial_{\bar \omega_h}$ for $h\geq n-l+1$ are non-singular at $M$. 
The key point is that the boundary integral in \eqref{2.10} is positive due to our assumption of $q$-pseudoconvexity as restated in Proposition~\ref{p1.1}. 
By discarding the positive boundary integrals we are thus reduced to integrals involving only terms of type $\partial_{\bar\omega_h}u_J\bar u_I$ for $h\geq n-l+1$. These latter are in turn reduced to terms of type $\partial_{\bar\omega_j}u_J\bar u_I$ for $j\leq n-l$ due to the choice of the distinguished representative of the form $u$ whose coefficients satisfy in particular \eqref{2.0}.   Summarizing up,                                                                                                                                                                                                \eqref{2.6} can be rewritten as
\begin{gather}
\Label{final1}
\sumK \left(\sumij\int_{U_\epsilon}e^{-\phi}\phi_{ij}u_{iK}\bar u_{jK} dV-\sumjq\int_{U_\epsilon}e^{-\phi}\phi_{jj}|u_J|^2 dV\right)
\\
\Label{final2}
+\sumK \left(\sumij\int_{\partial U_\epsilon}e^{-\phi}\rho_{ij}^\xi u_{iK}\bar u_{jK}-\sumjq\int_{\partial U_\epsilon}e^{-\phi}\rho_{jj}^\xi|u_J|^2 dV\right)
\\
\Label{final3}
+S+R^8,
\end{gather}
with  $R^8$ having the same estimate as prior error terms and with $\xi=\partial\tilde \rho$. Finally, by Proposition~\ref{p1.1}, the term in \eqref{final2} is bigger than $-O(\epsilon)\int_{\partial U_\epsilon}|u|^2 dV=-O(\epsilon^l)||u||^2_{H^0_\phi(M)}+o(\epsilon^l)=-c'||u||^2_{H^0_\phi(U_\epsilon)}+o(\epsilon^l)$. Note that the term in \eqref{final1} is bigger than $(t+c)||u||^2_\phi$ for large $t$. 
If we then choose $c$ which takes care of $c'$ and of the large constant for the estimate of $R^8$, we get from \eqref{2.2} the conclusion of the theorem.

\epf

\section{  Tangential estimates}
We recall that we are choosing an orthonormal basis of $(1,0)$ forms $\{\omega\}=\{\omega',\omega''\}$ satisfying
$$
\partial'_{\omega_j}\rho^h\equiv0,\,\,\partial''_{\bar\omega_j}\rho^h=\varkappa_{jh}.
$$
We recall that  $\partial'_{\omega_j}|_M$ and $\partial'_{\bar\omega_j}|_M$ for $j\leq n-l$ are the tangential vector fields of type $(1,0)$ and $(0,1)$ respectively and that the $\mathcal T_h:=\partial''_{\omega_h}-\partial''_{\bar\omega_h}$ and $\mathcal N_h:=\partial''_{\omega_h}+\partial''_{\bar\omega_h}$ for $h\geq n-l$ are the vector fields totally real tangential and normal to $M$ respectively. 
We also choose the extension of our forms $u$ from $M$ to $\C^n$ such that $u^\nu\equiv0$ and 
\begin{equation*}
\partial''_{\bar\omega_h}u_J=\sum_{j\leq n-l}a_j\partial_{\bar\omega_j}u_J+O(|\rho|),
\end{equation*}
for small coefficients $a_j$.
 By the $C^1$ regularity of the extensions, we then get
\begin{equation}
\begin{split}
\Label{3.1}
u^\tau &=u^\tau|_M+O(|\rho|),\,\,\partial'_{\omega_j}u^\tau=\partial'_{\omega_j}u^\tau|_M+
O(|\rho|),
\\
\mathcal T_ju^\tau &=\mathcal T_ju^\tau|_M+O(|\rho|). 
\end{split}
\end{equation}
We note that \eqref{3.1} implies for $u$ the follwing relations between its coefficients $u_J$ and their restrictions $(u_J)|_M$ 
\begin{equation}
\begin{split}
\Label{3.1bis}
||  u_J||_{H^0_\phi(U_\epsilon)} &=\epsilon^l||u_J||_{H^0_\phi(M)}+o(\epsilon^l),\\ 
||\partial_{\bar\omega_j}u_J||_{H^0_\phi(U_\epsilon)} &=\epsilon^l||\partial_{\bar\omega_j}u_J||_{H^0_\phi(M)}+o(\epsilon^l)
\end{split}
\end{equation}
 and so on. 
 We denote by $\bar\partial_b$ and $\bar\partial^*_b$ the tangential complexes to $M$ 
 associated to $\bar\partial$ and $\bar\partial^*$ respectively.  \eqref{3.1bis} immediately yields
\bl
\Label{l3.1}
In the above situation we have
\begin{gather}
\Label{3.2} ||u||_{H^0_\phi(U_\epsilon)}=\epsilon^l||u||_{H^0_\phi(M)}+o(\epsilon^l),
\\
\Label{3.3}   ||\bar\partial' u||_{H^0_\phi(U_\epsilon)}=\epsilon^l||\bar\partial_b u||^2_{H^0_\phi(M)}
+o(\epsilon^l),
\\
\Label{3.4}    ||\bar\partial^{\prime*}u||_{H^0_\phi(U_\epsilon)}=\epsilon^l||\bar\partial^*_bu||_{H^0_\phi(M)}+o(\epsilon^l).
\end{gather}
\el
\bpf
\eqref{3.2} is obvious. As for \eqref{3.3}, we have
$$
\bar\partial u=\sumJ{\underset{j\notin J}\sum}\partial_{\bar\omega_j}u_J\bar\omega_j\wedge\bar\omega_J,
$$
and
$$
\bar\partial_b u=\sumJ{\underset{j\notin J,\,j\leq n-q}\sum}\partial'_{\omega_j}u_J|_M\bar\omega'_j\wedge\bar\omega_J.
$$
Since
\begin{equation}
\partial'_{\omega_j}u_J=\partial'_{\omega_j}u_J|_M+o(|\rho|^l),
\end{equation}
then \eqref{3.3} immediately follows. Similarly
$$
\bar\partial^*u=-\sumK\sumj\delta_{\omega_j}u_{jK}\bar\omega_K,
$$
and
$$
\bar\partial^*_bu=\sumK\underset {j\leq n-l}\sum \delta'_{\omega_j}u_{jK}\bar\omega_K,
$$
where we remember that $\delta_{\omega_j}=\partial_{\omega_j}-\phi_j$.
(Note here that $\bar\partial^*_b=\bar\partial^*|_M$ over $\bar\partial$-Neumann forms.)
\epf 
We go back to Theorem~\ref{t2.1}. 
We recall that $U_\epsilon$ denotes the intersection of the tube defined by $\tilde\rho<0$ with a ball $B$ centered at $z_o$; we will consider the neighborhood of $z_o$ defined by $M'=M\cap B$. 
If we multiply both sides of \eqref{2.1} by $\epsilon^{-l}$ and go to the limit for $\epsilon\to0$ we get for any large $t$ and for any tangential form $u$ of degree $k\geq q+1$
\begin{equation}
\Label{tangential}
\frac t3||u||^2_{H^0_\phi(M')}\leq ||\bar\partial_bu||^2_{H^0_\phi(M')}+||\bar\partial^*_bu||^2_{H^0_\phi(M')}.
\end{equation}
We deal  now with the (unweighted) Sobolev spaces $H^s$ (for $s$ integer). We will emphasize from now on the dependence of $\phi_t$ on $t$. 
We will assume also  that $M$ is $C^\infty$.
The main result of the section is the following
\bt
\Label{t3.1}
Let $M$ be $q$-pseudoconvex at $z_o$. Then
for any $s$, for any sufficiently large $t= t_s$,   for suitable $c=c_{t_s}$ and for a suitable neighborhood $M'$ of $z_o$  
we have
\begin{equation}
\Label{sobolev0}
||u||^2_{H^s(M')}\leq c(||\bar\partial_bu||^2_{H^s(M')}+||\bar\partial^*_{b,t_s}u||^2_{H^s(M')}),
\end{equation}
for any tangential form $u$ of length $k\geq q+1$, resp. $k\leq p-1$. (Here we write $\bar\partial^*_{b,\phi_{t_s}}$ to emphasize the dependence on the weight $\phi_{t_s}$.)
\et
Note that the weight $\phi_{t_s}$, which is eliminated in the norms, reappears in an essential way in the operation of adjunction.
\bpf
We denote by a common symbol $\mathcal T$ all tangent vector fields that is any combination of the $\partial'_{\omega_j}$'s, $\partial'_{\bar\omega_h}$'s and $(\partial''_{\omega_h}-\partial''_{\bar\omega_h})$'s. 
If $\alpha$ is a multiindex, we set $\mathcal T^\alpha=\mathcal T_1^{\alpha_1}\dots\mathcal T_{n-l}^{\alpha_{n-l}}$. We write the commutators $[\bar\partial^*_{b,\phi_{t_s}},\mathcal T^\alpha]=A_s+A^t_{s-1}$ where $A_s$ is an operator of degree $s$ independent of $t$ and $A^t_{s-1}$ is of degree $s-1$; thus the  coefficients of $A^t_{s-1}$ are estimated by $t$. It follows that $||A_su||_{H^0}\leq a_s||u||_{H^s}$ and $||A^t_{s-1}u||_{H^0}\leq a_s t||u||_{H^{s-1}}$ for a suitable constant $a_s$. 
We apply \eqref{tangential} to all terms of the type $\mathcal T^\alpha u$ for $|\alpha|=s$; we have
\begin{equation}
\Label{sobolev1}
\begin{split}
\frac t3||\mathcal T^\alpha u||^2_{H^0}&\leq ||\bar\partial_b\mathcal T^\alpha u||^2_{H^0}+||\bar\partial^*_b\mathcal T^\alpha u||^2_{H^0}
\\
&\leq ||\mathcal T^\alpha\bar\partial_bu||^2_{H^0}+||\mathcal T^\alpha\bar\partial^*_bu||^2_{H^0}+a_s||u||_{H^s}+a_st||u||^2_{H^{s-1}}.
\end{split}
\end{equation}
Now, by inductive assumption we have 
$$
a_s||u||^2_{H^s}+a_st||u||^2_{H^{s-1}}\leq a_s||u||^2_{H^s}+
a_sc_{s-1}t(||\bar\partial_bu||^2_{H^{s-1}}+||\bar\partial^*_bu||^2_{H^{s-1}}).
$$
If we take $t$ so large that $\frac t3-a_sc_{s-1}\geq1$ (in such a way that the term involving $||u||_{H^s}$ in the right side of \eqref{sobolev1} can be ``absorbed" in the left), and define $c_s:=1+a_sc_{s-1}t$,  we get the conclusion of the proof of \eqref{sobolev0}.

\epf
We point out that only the use of the weight $\phi_{t_s}$ produces a big constant on the left side of \eqref{tangential} which makes it possible to pass through derivatives absorbing the constants $a_s$
and $c_{s-1}$ in the above proof. Once this is carried out, we come back to unweighted estimates (since the spaces $H^0$ and $H^0_{\phi_{t_s}}$ coincide and have equivalent norms). Thus, we did eventually got rid of the weights from our norms. However they did a great service and gave the control of the derivatives of the coefficients of our forms $u$. 

%
%
%%%%%%%%%%%%%%%%%%%%%%%%%%%%%%%%%%%%%%%%%%%%%%%%%%%%%%%%%%%%%%%%%%%%%%%%%%%%%%%%%%%%
\section{Existence theorems for $\bar\partial_b$}
The main applications of the tangential estimates of \S~3 consist in results of local existence of $C^\infty$ solutions for $\bar\partial_b$. We will follow here closely the theory by {\sc Kohn}.
If $s$ is any Sobolev index, we take $t=t_{s}$ such that the conclusions of Theorem~\ref{t3.1} hold: thus \eqref{sobolev0} is satisfied. 
We recall that we are denoting by $\bar\partial^*_{b,t_s}$   the transposed of $\bar\partial_b$ in the $H^0_{\phi_{t_{s}}}$  scalar product. We set
$$
\Box_{b,t_s}=\bar\partial_b\bar\partial^*_{b,t_s}+\bar\partial^*_{b,t_s}\bar\partial_b.
$$
We remark that with this notation, \eqref{sobolev0} can be rewritten as
\begin{equation}
\Label{continuousbox}
\begin{split}
\frac t3||u||^2_{H^0_{\phi_{t_s}}}&\leq(\Box_{b,t_s}u,u)_{H^0_{\phi_{t_s}}}
\\
&\leq (\Box_{b,t_s}u,u)_{H^0_{\phi_{t_s}}}\leq||\Box_{b,\phi_{t_s}}u||_{H^0_{\phi_{t_s}}}||u||_{H^0_{\phi_{t_s}}},
\end{split}
\end{equation}
for any tangential form $u$ of degree $k\geq q+1$. 
Denote by $\T{R}_{\Box_{b,t_s}}$ and $D_{\Box_{b,t_s}}$ the range and the domain of $\Box_{b,t_s}$ respectively.
It follows from \eqref{continuousbox} that $\T{R}_{\Box_{b,t_s}}$ is closed and $\Box_{b,t_s}$ is injective. From the orthogonal decomposition $H^0_{\phi_{t_s}}=\T{R}_{\Box_{b,t_s}}\oplus \T{Ker}\,\Box_{b,t_s}=\T{R}_{\Box_{b,t_s}}$, we conclude that there is a well defined ``weighted $\bar\partial$-Neumann operator" 
$$
N_{b,t_s}\,:\,L^2\to D_{\Box_{b,t_s}},
$$
such that $N_{b,t_s}\Box_{b,t_s}=\Box_{b,t_s}N_{b,t_s}=\T{id}$ and which satisfies  
\begin{equation}
\Label{neumanncontinuous}
t||N_{b,t_s}f||_{H^0_{\phi_{t_s}}}\simleq||f||^2_{H^0_{\phi_{t_s}}}\quad\forall f\in C^\infty,
\end{equation}
where ``$\simleq$" denotes estimation up to a multiplicative constant. 
We can also rephrase the conclusions of Theorem~\ref{t3.1} in terms of the weighted Neumann operator: for any $s$ and for a suitable $t_s$ we have
\begin{equation}
\Label{apriori}
||N_{b,t_s}f||_{H^s}\simleq||f||_{H^s}\T{ if $f$ and $N_{b,t_s}f$ are $C^\infty$}.
\end{equation}
We want to get rid of the condition $N_{b,t_s}f\in C^\infty$ from equation \eqref{apriori}. 
For this purpose we define  an elliptic perturbation $\Box^\sigma_{b,t_s}:=\Box_{b,t_s}+\sigma(\sum\mathcal T^2)$ 
where the sum is extended to a full set of tangential vector fields. 
This yields an inverse ``regularizing" operator
\begin{equation}
\Label{4.2}
N^\sigma_{b,t_s}\,:\,H^s\to D_{\Box_{b,t_s}}\cap H^{s+1},
\end{equation}
which satisfies
\begin{equation}
\Label{4.3}
||N^\sigma_{b,t_s}f||_{H^s}+\sigma||N^\sigma_{b,t_s}f||_{H^{s+1}}\simleq||f||_{H^s}.
\end{equation}
It follows that for some  $\sigma_j\to0$, the sequence $N^{\sigma_j}_{b,t_s}f$ has a weak $H^s$-limit. Hence for $f\in H^s$, we have  $N_{b,t_s}f\in H^s$ and $N^{\sigma_j}_{b,t_s}\to N_{b,t_s}f$; in particular
\begin{equation}
||N_{b,t_s}f||_{H^s}\simleq||f||_{H^s}\quad\forall f\in C^\infty,
\end{equation}
and, in fact, for any $f\in H^s$ by density.
By using the above construction we get the following statement
\bp
\Label{p4.1}
Fix $s$ and assume that for suitable $t=t_{s}$ \eqref{sobolev0} holds 
for forms of a certain degree $k$. 
Let $f$ be a $C^\infty$ form on $M'$ of degree $k$ satisfying $\bar\partial_b f=0,$ and define $u:=\bar\partial^*_bN_{b,t_s}f$; then $u$ belongs to $H^s$ and satisfies
\begin{equation}
\Label{4.4}
\begin{cases}
\bar\partial_bu=f,
\\
||u||_{H^s}\simleq ||f||_{H^s}.
\end{cases}
\end{equation}

\ep
The afore-defined $u$ is orthogonal to $\T{Ker }\bar\partial_b$ and it is also clear 
that under such condition there is uniqueness for the solution.
Note  that according to Theorem~\ref{t3.1}, the hypotheses of Proposition~\ref{p4.1} are 
fulfilled, for any $s$, for any degree $k\geq q+1$ and for a suitable neighborhood $M'$ of $z_o$,  when $M$ is 
$q$-pseudoconvex at $z_o$.
\bpf

We have
\begin{equation}
\Label{4.5}
\begin{split}
\bar\partial_bg &=\bar\partial_b\Box_{b,t_s}N_{b,t_s}g=\bar\partial_b 
\left(\bar\partial^{*}_{b,t_s}\bar\partial_b+\bar\partial_b\bar\partial^{*}_{b,t_s}\right)
N_{b,t_s}g
\\
&=\left(\bar\partial_b\bar\partial^{*}_{b,t_s}+\bar\partial^{*}_{b,t_s}\bar\partial_b\right)
\bar\partial_bN_{b,t_s}g.
\end{split}
\end{equation}
Hence, if $f$ satisfies $\bar\partial_bf=0$, we have
\begin{equation}
\begin{split}
0=N_{b,t_s}\bar\partial_bf&=N_{b,t_s}\Box_{b,t_s}\bar\partial_bN_{b,t_s}f
\\
&=\bar\partial_bN_{b,t_s}f.
\end{split}
\end{equation}
It follows that for $u:=\bar\partial^{b*}N_{b,t_s}f$ we have
\begin{equation*}
\begin{split}
f&=\left(\bar\partial_b\bar\partial^{*}_{b,t_s}+\bar\partial^{*}_{b,t_s}\bar\partial_b\right)N_{b,t_s}f
\\
&=\bar\partial_b(\bar\partial^{*}_{b,t_s}N_{b,t_s}f)
\\
&=\bar\partial_b(\bar\partial^{*}_{b,t_s}N_{b,t_s}f).
\end{split}
\end{equation*}
This completes the proof of the first of \eqref{4.4}.

As for the second, we first recall that $||N^b_{s,t_s}f||_{H^s_{\phi_{t_s}}}\simleq||f||_{H^s_{\phi_{t_s}}}$. Next, we remark that
\begin{gather}
(\bar\partial_bN_{b,t_s}f,\bar\partial_bN_{b,t_s}f)+(\bar\partial^{*}_{b,t_s}N_{b,t_s}f,\bar\partial^{*}_{b,t_s}N_{b,t_s}f)
\\
=(\Box_{b,t_s}N_{b,t_s}f,N_{b,t_s}f)
=(f,N_{b,t_s}f)
\\
\leq||f||_{H^s_{\phi_{t_s}}}||N_{b,t_s}f||_{H^s_{\phi_{t_s}}}\simleq ||f||^2_{H^s_{\phi_{t_s}}}.
\end{gather}
This implies
immediately the second of  \eqref{4.4}.
\epf
Once we know that for $f$ in $C^\infty$ with $\bar\partial_bf=0$ we can find a solution of $\bar\partial_bu=f$ with estimate in each $H^s$, it is easy to see that we can find indeed a  solution in $C^\infty$. The proof consists in a variant of a very classical approximation argument due to H\"ormander (as referred by Kohn in \cite{K73}).
\bt
\Label{t4.1}
Let $M$ be $C^\infty$ and $q$-pseudoconvex at $z_o$. Then for any $f$ in $C^\infty$ of degree 
$k\geq q+1$ with $\bar\partial_bf=0$, we can find a $C^\infty$ solution $u$   of 
$\bar\partial_bu=f$ at $z_o$.
\et
\bpf
According to Theorem~\ref{t3.1}, for any $s$ and for suitable $t=t_{s}$ \eqref{sobolev0} holds 
for forms of any degree $k\geq q+1$
and for  a suitable neighborhood $M'$ of $z_o$. 
 According to Proposition~\ref{p4.1} we can find for any $s$ 
an $H^s$ solution $u_s$ in $M'$ with the estimate \eqref{4.4}. We want to carry on our proof by showing 
 by induction that there is a sequence of solutions $u_\nu\in H^\nu$ of 
$\bar\partial_bu_\nu=f$ which satisfies
\begin{equation}
\Label{4.6}
||u_{\nu+1}-u_\nu||_{H^{\nu}}\leq2^{-\nu}.
\end{equation}
In fact, once $u_1,\dots,u_\nu$ have been found, we take $\tilde u_{\nu+1}\in C^\infty$ and 
$v_{\nu+1}\in H^{\nu+1}$ such that
\begin{equation}
\Label{4.7}
\begin{cases}
||\tilde u_{\nu+1}-u_\nu||_{H^\nu}\leq2^{-(\nu+1)}\T{ and }||\bar\partial_b\tilde u_{\nu+1}-\bar\partial_bu_\nu||_{H^\nu}<2^{-(\nu+1)},
\\
\bar\partial_bv_{\nu+1}=f-\bar\partial_b\tilde u_{\nu+1},
\\
\begin{split}
||v_{\nu+1}||_{H^{\nu}}&\leq||f-\bar\partial_b\tilde u_{\nu+1}||_{H^{\nu}}
\\
&=||\bar\partial_bu_\nu-\bar\partial_b\tilde u_{\nu+1}||_{H^{\nu}}
\\
&\leq2^{-(\nu+1)}.
\end{split}
\end{cases}
\end{equation}
If we then set 
$$
u_{\nu+1}:=\tilde u_{\nu+1}+v_{\nu+1},
$$
we have
\begin{equation}
\begin{split}
||u_{\nu+1}-u_\nu||_{H^{\nu}}&\leq||\tilde u_{\nu+1}-u_\nu||_{H^{\nu}}+||v_{\nu+1}||_{H^{\nu}}
\\
&\leq2^{-(\nu+1)}+2^{-(\nu+1)}=2^{-\nu}.
\end{split}
\end{equation}
Thus $u_{\nu+1}$ is a solution of $\bar\partial_bu_{\nu+1}=f$ in $M'$ which satisfies  \eqref{4.6}. 
\epf

\section{The hypersurface case}
We want to end by discussing in greater detail the case of a hypersurface $M$. We have already seen that by setting $q=\T{max}(s^-(z_o)+s^0(z_o),s^+(z_o)+s^0(z_o))$ we have local $q$-pseudoconvexity at $z_o$. Also, if $s^0$ is locally constant when we move $z$, then it has been proved that in fact local $q$-pseudoconvexity holds for the lower choice $q=\T{max}(s^-(z_o),s^+(z_o))$. In both cases the equation $\bar\partial_bu=f$ has local $C^\infty $ solution $u$ for any $C^\infty$ datum $f$ with $\bar\partial_b f=0$ in any degree $k$ bigger than the corresponding $q$. However, we can improve much our existence theorems. To this end we denote by $U^\pm$ the two components of $\C^n\setminus M$ with outward conormals $\pm\xi$. We still assume that $s^\pm$ are constant and notice that $s^-(z,-\xi)=s^+(z,+\xi)$. The argument of the above sections can be applied separately to each domain $\bar U^+$ and $\bar U^-$ which is $s^-$ and $s^+$ pseudoconvex with respect to its respective conormal:
\bp
\Label{p5.1}
Let $s^0$ be constant in a neighborhood of $z_o$. Then for any $f$ with  $C^\infty(\bar U^\pm)$ coefficients in
a  neighbourhood of $z_o$,
 satisfying $\bar \partial_bf=0$ and of degree $k\geq s^\mp+1$, there exists a solution $u$ of $\bar\partial_bu=f$ in a neighbourhood of $z_o$ with coefficients in $C^\infty(\bar U^\pm)$.
\ep
We pass to consider the equation $\bar\partial_bu=f$  for $k\leq s^\mp-1$.  In this case, we have the so called local $s^\mp$-pseudoconcavity. Similar arguments as in Section 2 go through without need of a weight $t|z|^2$ and yield the so called ``subelliptic estimates"
\begin{multline}
\Label{5.1}
||u||^2_{H^{\frac12}(\bar U^\pm)}\underset\sim< ||\bar\partial u||^2_{H^0(\bar U^\pm)}+||\bar\partial^* u||^2_{H^0(\bar U^\pm)}+||u||^2_{H^0(\bar U^\pm)}
\\ \forall u
\T{ of degree $k\leq s^\mp-1$.}
\end{multline}
The estimate \eqref{5.1} yields ``gain" of regularity for the solution $u$ with respect to the datum $f$; in particular it implies the hypoellipticity of the system $(\bar\partial,\,\bar\partial^*)$. Moreover, by replacing in the calculations of Section 2 the weight $t|z|^2$  by $-t\sum_{j\leq s^-}|z_j|^2+t\sum_{j\geq s^-+1}|z_j|^2$ in case of $\bar U^+$, resp. $-t\sum_{j\leq s^+}|z_j|^2+t\sum_{j\geq s^++1}|z_j|^2$ for $\bar U^-$, we can prove $H^0$ estimates for $\bar U^+$, resp. $\bar U^-$, of the type of those in Theorem~\ref{t2.1}, which imply local existence of $H^0$ solutions (cf. \cite{H65} Theorem 3.3.1). In combination with the afore-mentioned hypoellipticity this implies that the equation $\bar\partial u=f$ with $\bar\partial f=0$ is locally solvable in $C^\infty(\bar U^\pm)$ for any degree $k\leq s\mp-1$.
On the other hand it is classical that the tangential $\bar\partial$-problem for the hypersurface $M$ can be split into the $\bar\partial$ problems for the half-spaces $\bar U^\pm$. In fact, any germ of $C^\infty$ form $f$ satisfying $\bar\partial_bf=0$ on $M$ can be decomposed into the sum $f=f^+\oplus f^-$ with $f^\pm$ satisfying $\bar\partial f^\pm=0$ on $\bar U^\pm$. This yields
\bt
\Label{t5.1}
(Cf. \cite{Z02})
Let $M$ be a $C^\infty$ hypersurface such that $s^0$ is constant in a neighborhood of $z_o$. Then for any germ of $C^\infty(M)$ form $f$ at $z_o$ of degree $k\neq s^-,\,\,s^+$, satisfying $\bar\partial_bf=0$, there exists a germ of $C^\infty(M)$ form $u$ which solves $\bar\partial_bu=f$.
\et
Let us point out that according to \cite{Z02}  the equation $\bar\partial_bu=f$ is not solvable in the two critical degrees $s^-$ and $s^+$;
 when the Levi form of $M$ is non-degenerate, that is $s^0=0$, the result was alredy proved in \cite{AFN81}. 
If we go back to the literature, the solvability of the system $\bar\partial_b$ in degree $k$ is related to the so-called $Y(k)$-condition by {\sc Kohn} and {\sc H\"ormander}:
the Levi form of $M$ has $\T{max}(k+1,n-k)$ eigenvalues of the same sign or $\T{min}(k+1,n-k)$ pairs of eigenvalues of opposite sign at each point. 
Another equivalent formulation of $Y(k)$ is that: $k\notin[s^-,s^-+s^0]\cup[s^+,s^++s^0]$. 
Under this condition Kohn and H\"ormander proved tangential estimates of subelliptic type \eqref{5.1} which yield existence of  smooth solutions in degree $k$, except for a finite-dimensional set of $f$. Indeed, by an argument similar to the one which led to Theorem~\ref{t5.1}, they proved that there are no exceptions at all. If we compare with our Theorem~\ref{t5.1}, we see that when $s^0$ is constant, then we have got new results of solvability for all indices $k\in(s^-,s^-+s^0]\cup(s^+,s^++s^0]$. 
 
%\epf


\begin{thebibliography}{AG62}
\bibitem[1]{A03} H. Ahn, Global boundary regularity of the $\bar\partial$-equation on
$q$-pseudoconvex domains, {\em Preprint} (2003)
%\bibitem[2]{AG62} A. Andreotti, H. Grauert, Th\'eor\`emes de finitude 
%pour la cohomologie des \'espaces complexes,
%{\em Bull. Soc. Math. de France}
%{\bf 90} (1962), 193--259
%
\bibitem[2]{AFN81} A. Andreotti, G. Fredricks, M. Nacinovich, On the absence of a Poincar\'e lemma in tangential Cauchy-riemann complexes, {\em Ann. Scuola Norm. Sup. Pisa} {\bf 8} (1981), 365-404

\bibitem[3]{BZ03} L. Baracco, G. Zampieri, Global regularity for 
$\bar\partial$ on $q$-pseudoconvex domains, {\em Preprint} (2003)
%
\bibitem[4]{B92} D. Barrett, Behavior of the Bergman projection on 
the Diederich-Fornaess worm, {\em Acta Math.} {\bf 168} (1992), 1--10
%
\bibitem[5]{CS01} S.C. Chen, M.C. Shaw, Partial differential equations in several complex variables, {\em Studies in Adv. Math. - AMS Int. Press} {\bf 19} (2001)

\bibitem[6]{Ch96} M. Christ, Global $C^\infty$ irregularity of the 
$\bar\partial$-Neumann problem for worm domains, {\em J. of the 
A.M.S.}{\bf 9}-(4) (1996), 1171--1185
\bibitem[7]{D78} M. Derridj, Regularit\'e pour $\bar\partial$ dans quelques domaines faiblement pseudo-convexes, {\em J. Differential Geometry} {\bf 13} (1978), 559-576
\bibitem[8]{DT76} M. Derridj, D. Tartakoff, Sur la r\'egularit\'e 
analytique globale des solutions du probl\`eme de Neumann pour 
$\bar\partial$, {\em S\'em. Goulaouic-Schwartz}, (1976)

%
\bibitem[9]{Du79} A. Dufresnoy, Sur l'operateur $\bar\partial$ et les 
fonctions difff\'erentiables au sens  de Whitney, {\em Ann. de 
l'Inst. Fourier} {\bf (29)}-(1) (1979), 229--238
\bibitem[10]{Hen77} Henkin G.M., H. Lewy's equation and analysis on 
pseudoconvex manifolds (Russian), I, {\em Uspehi Mat. Nauk.} {\bf 
32}-(3) (1977), 57--118
%
\bibitem[11]{Ho92} L.H. Ho, $\bar\partial$-problem on weakly $q$-convex domains, {\em Math. Ann.} {\bf 290} (1) (1991), 3--18
%
\bibitem[12]{H65} L. H\"ormander, $L^2$ estimates and existence 
theorems for the $\bar\partial$ operator,
{\em Acta Math.} {\bf 113} (1965), 89-152
\bibitem[13]{H73} L. H\"ormander, An introduction to complex analysis 
in several complex variables, Van Nostrand, Princeton N.J., 1966
\bibitem[14]{K73} J.J. Kohn, Global regularity for $\bar\partial$ on 
weakly pseudo-convex manifolds, {em Transactions of the A.M.S.} {\bf 
181} (1973), 273--292
\bibitem[15]{K77} J.J. Kohn, Methods of partial differential 
equations in complex analysis,
{\em Proceedings of Symposia in pure Mathematics} {\bf 30} 
(1977),215--237
%
\bibitem[16]{K79} J.J. Kohn, Subellipticity of the $\bar\partial$-Neumann problem on pseudoconvex domains: sufficient conditions, {\em Acta Math.} {\bf 142} (1979), 79--122


%\bibitem[17]{MS98} J. Michel, M.C. Shaw, $C^\infty$ regularity of solutions of the tangential CR-equations on weakly %pseudoconvex manifolds, {\em Mat. Ann.} {\bf 311} (1998), 147-162.

\bibitem[17]{N84} M. Nacinovich, Poincar\'e lemma for the tangential Cauchy Riemann complexes, {\em Math. Ann.} {\bf 268} (1984), 449-471

\bibitem[18]{S92} M.C. Shaw, Local existence theorems with estimates for $\bar\partial_b$ on weakly pseudoconvex boundaries, {\em Mat. Ann.} {\bf 294} (1992), 677-700 
\bibitem[19]{Z00} G. Zampieri, $q$-Pseudoconvexity and regularity at 
the boundary for solutions of the $\bar\partial$-problem, {\em 
Compositio Math.} {\bf 121} (2000), 155--162
\bibitem[20]{ZII00} G. Zampieri, Solvability of the $\bar\partial$ 
problem with $C^\infty$ regularity up to the boundary on wedges of 
$\C^N$,
{\em Israel J. of Math.} {\bf 115} (2000), 321--331
\bibitem[21]{Z02} G. Zampieri, $q$=pseudoconvex hypersurfaces through higher codimensional submanifolds of $\C^n$ {\em J. Reine Angew Math.} {\bf 544}, 83-90
\end{thebibliography}
\end{document}